\documentclass[reqno,draft]{amsart}

\theoremstyle{plain}
\newtheorem{theorem}{Theorem}
\newtheorem{lemma}[theorem]{Lemma}

\newtheorem{proposition}[theorem]{Proposition}

\theoremstyle{remark}

\renewcommand{\theequation}{\arabic{section}.\arabic{equation}}
\renewcommand{\thetheorem}{\arabic{section}.\arabic{theorem}}

\def\aa{\alpha}

\def\Dd{\Delta}

\def\pp{\partial}
\def\na{\nabla}

\begin{document}
\title[Two regularity criteria for the 3D MHD equations]{Two regularity criteria for the 3D MHD equations}

\author{Chongsheng Cao}
\address{Department of Mathematics, Florida International University, Miami, FL 33199, USA}
\curraddr{}
\email{caoc@fiu.edu}
\thanks{}

\author{Jiahong Wu}
\address{Department of Mathematics, Oklahoma State University, Stillwater, OK 74078, USA}
\curraddr{}
\email{jiahong@math.okstate.edu}
\thanks{}



\date{}

\dedicatory{}

\begin{abstract}
This work establishes two regularity criteria for the 3D incompressible MHD equations. The first one is in terms of the derivative of the velocity field in one-direction while the second one requires suitable boundedness of the derivative of the pressure in one-direction.
\end{abstract}

\maketitle

\vskip 0.1in
\section{Introduction}

\setcounter{section}{1}
\setcounter{theorem}{0}

This paper is concerned with the global regularity of solutions to the 3D incompressible magneto-hydrodynamical (MHD) equations
\begin{eqnarray}
&&\hskip-.8in
u_t + u\cdot\nabla u =\nu\, \Delta u -\nabla p +  b\cdot\nabla b, \qquad x\in {\mathbf R}^3, \; t >0,  \label{DE-1}  \\
&&\hskip-.8in
b_t + u\cdot\nabla b = \eta\, \Dd b + b\cdot\nabla u, \qquad x\in {\mathbf R}^3, \; t >0, \label{DE-2}   \\
&&\hskip-.8in
\nabla\cdot u=0,  \qquad x\in {\mathbf R}^3, \; t >0, \label{DE-3} \\
&&\hskip-.8in
\nabla \cdot b=0, \qquad x\in {\mathbf R}^3, \; t >0, \, \label{DE-4}
\end{eqnarray}
where $u$ is the fluid velocity, $b$ the magnetic field, $p$ the pressure, $\nu$ the viscosity and $\eta$ the magnetic diffusivity. Without loss of generality, we set $\nu=\eta=1$ in the rest of the paper. The MHD equations govern the dynamics of the velocity and magnetic fields in electrically conducting fluids such as plasmas. (\ref{DE-1}) reflects the conservation of momentum, (\ref{DE-2}) is the induction equation and (\ref{DE-3}) specifies the conservation of mass. Besides their physical applications, the MHD equations are also mathematically significant. Fundamental mathematical issues such as the global regularity of their solutions have generated extensive research and many interesting results have been obtained (see, e.g., \cite{RuSc},\cite{CKS},\cite{CaoWu},\cite{Chae},\cite{Chan},\cite{CMZ},
\cite{CorMar},\cite{DL},\cite{GiOh},\cite{HeWang},\cite{HeXin1},\cite{HeXin2},\cite{LeiZhou},
\cite{MiYu},{\cite{MiYuZh},\cite{NM},\cite{Oh},\cite{PoPoSu},\cite{ScScSu},\cite{SeTe},
\cite{Wu1},\cite{Wu2},\cite{Wu3},\cite{Wu4},\cite{Wu5},\cite{Yuan},\cite{Zhou2}).

\vskip 0.1in
Attention here is focused on the global regularity of solutions to the initial-value problem (IVP) of (\ref{DE-1}),(\ref{DE-2}),(\ref{DE-3}) and (\ref{DE-4}) with a given initial data
\begin{eqnarray}
&&\hskip-.8in
u(x,0)=u_0(x), \quad b(x,0)=b_0(x), \qquad x\in {\mathbf R}^3.
\end{eqnarray}
It is currently unknown whether solutions of this IVP can develop finite time singularities even if $(u_0, b_0)$ is sufficiently smooth. This work presents new regularity criteria under which the regularity of the solution is preserved for all time. The global regularity issue has been thoroughly investigated for the 3D Navier-Stokes equations and many important regularity criteria have been established (see, e.g., \cite{BKM},\cite{Beirao1},\cite{Beirao2},\cite{BeGa},\cite{Cao},\cite{CQT},\cite{CaoT},\cite{Chae2},
\cite{Chae3},\cite{ChaeCh},\cite{ESS1},\cite{ESS2},\cite{He1},\cite{KoOg},\cite{KoYa},
\cite{KuMo1},\cite{KuMo2},\cite{NePe},\cite{PePo},\cite{Po},\cite{SeSv},\cite{Serr},\cite{ZhCh},\cite{Zhou1}).  Some of these criteria can be extended to the 3D MHD equations by making assumptions on both $u$ and $b$ (see, e.g., \cite{CKS},\cite{Wu1}). Realizing the dominant role played by the velocity field in the regularity issue, He and Xin was able to derive criteria in terms of the velocity field $u$ alone (\cite{HeXin1},\cite{HeXin2}). They showed that, if $u$ satisfies
\begin{eqnarray}
&&\hskip-.8in
\int_0^T \|\nabla u(\cdot, t)\|_\alpha^\beta \,dt <\infty \quad\mbox{with}\quad
\frac{3}{\alpha}+ \frac{2}{\beta} =2\quad\mbox{and}\quad 1<\beta \le 2,
\end{eqnarray}
then the solution $(u,b)$ is regular on $[0,T]$. This assumption was weakened in \cite{Wu5} with $L^\alpha$-norm replaced by norms in Besov spaces and further improved by Chen, Miao and Zhang in \cite{CMZ}. As pointed out in \cite{HeXin1},
the regularity criteria in terms of the velocity field alone are consistent with the numerical simulations in \cite {PoPoSu} and with the observations of space and laboratory plasmas in \cite{Has}.

\vskip .1in
This paper presents two new regularity criteria. The first one assumes
\begin{eqnarray}
&&\hskip-.8in
\int_0^T \|u_z(\cdot, t)\|_\alpha^\beta \,dt <\infty \quad\mbox{with}\quad \alpha\ge 3
\quad\mbox{and}\quad \frac{3}{\alpha}+ \frac{2}{\beta} \le 1  \label{vs}
\end{eqnarray}
and the second requires the pressure satisfy
\begin{eqnarray}
&&\hskip-.8in
\int_0^T \|p_z(\tau)\|^{\beta}_{\aa} \,d\tau < \infty \quad \mbox{with}\quad \aa \geq \frac{12}{7}\quad\mbox{and}\quad  \frac{3}{\aa} + \frac{2}{\beta} \leq \frac{7}{4}.
\label{ps}
\end{eqnarray}
That is, any solution $(u,b)$ of the 3D MHD equations is regular if the derivative of $u$ in one direction, say along the $z$-axis, is bounded in $L^\beta([0,T]; L^\alpha)$ with $(\alpha, \beta)$ satisfying (\ref{vs}) or if  the derivative of $p$ in one direction satisfies (\ref{ps}). The proof of the first criterion is accomplished through two stages with the first controlling the time integrals of $\|\na u_z\|_2$ and $\|\na b_z\|$ in terms of the $L^\beta([0,T]; L^\alpha)$-norm of $u_z$ and the second bounding $\|\na u\|_2$ and $\|\na b\|_2$ by the time integrals of $\|\na u_z\|_2$ and $\|\na b_z\|$. The details are presented in the second section. The criterion in terms of $p_z$ and its proof are provided in the third section.

\vskip .1in
We will use the following elementary inequalities:
\begin{equation}\label{3q}
\|\phi\|_{\gamma}  \le C\; \|\phi_x\|_\lambda^\frac13\; \|\phi_y\|_\lambda^\frac13 \; \|\phi_z\|^{\frac13}_\mu,
\end{equation}
where the parameters $\mu$, $\lambda$ and $\gamma$ satisfy
$$
1\le \mu, \;\lambda  <\infty, \quad 1<\frac{1}{\mu}+ \frac{2}{\lambda} \le 4 \quad
\mbox{and}\quad \gamma = \frac{3\lambda}{2-\lambda\left(1-\frac{1}{\mu}\right)},
$$
and
\begin{equation}\label{SL}
\|\phi\|_r  \le  C(r) \;\|\phi\|_2^{\frac{6-r}{2r}}\;  \|\phi_x\|_2^{\frac{r-2}{2r}} \;  \|\phi_y\|_2^{\frac{r-2}{2r}} \;  \|\phi_z\|_2^{\frac{r-2}{2r}}, \quad 2\le r\le 6.
\end{equation}
These inequalities may be found in \cite{Ad},\cite{Gal},\cite{Lad}. For the convenience of the readers, the proofs of these inequalities are provided in Appendix A. Throughout the rest of this paper the $L^p$-norm of a function $f$ is denoted by $\|f\|_p$, the $H^s$-norm by $\|f\|_{H^s}$ and the norm in the Sobolev space $W^{s,p}$ by $\|f\|_{W^{s,p}}$.

\vskip 0.3in
\section{Criterion in terms of $u_z$}

\setcounter{equation}{0}
\setcounter{theorem}{0}

This section establishes the regularity criteria in terms of $u_z$.

\begin{theorem}
Assume $(u_0, b_0)\in H^3$, $\na \cdot u_0=0$ and $\na\cdot b_0=0$. Let $(u,b)$ be the corresponding solution of the 3D MHD equations (\ref{DE-1}),(\ref{DE-2}),(\ref{DE-3}) and (\ref{DE-4}).  If $u$ satisfies
\begin{equation}\label{uzcon}
M(T) \equiv \int_0^T \|u_z(\cdot, t)\|_{\alpha}^\beta \;dt <\infty \quad\mbox{with}\quad \alpha \ge 3
\quad\mbox{and}\quad \frac{3}{\alpha}+\frac2\beta\le 1
\end{equation}
for some $T>0$, then $(u,b)$ can be extended to the time interval $[0,T+\epsilon)$ for some $\epsilon>0$.
\end{theorem}

The proof of this theorem is divided into two major parts. The first part establishes bounds for $\|u_z\|_2$, $\|b_z\|_2$ and the time integrals of $\|\na u_z\|_2^2$ and $\|\na b_z\|_2^2$ while the second controls  $\|\na u\|_2$ and $\|\na b\|_2$ in terms of the time integrals of $\|\na u_z\|^2_2$ and $\|\na b_z\|^2_2$.

\subsection{Bounds for $\|u_z\|_2$ and $\|b_z\|_2$} This subsection bounds $\|u_z\|_2$ and $\|b_z\|_2$ in terms of $M$ in (\ref{uzcon}).
\begin{proposition} \label{ma1}
Assume $(u_0, b_0)\in H^3$, $\na \cdot u_0=0$ and $\na\cdot b_0=0$. Let $(u,b)$ be the corresponding solution of the 3D MHD equations (\ref{DE-1}),(\ref{DE-2}),(\ref{DE-3}) and (\ref{DE-4}). Suppose (\ref{uzcon}) holds. Then, for any $t\le T$,
\begin{eqnarray}
&& \|u_z(t)\|_2^2+ \|b_z(t)\|_2^2 \le C\; e^{(\|u_0\|_2^2 +\|b_0\|_2^2)} \;e^{M(t)}\;
\nonumber\\
&& \qquad \qquad \times\left[\left(\|u_z(0)\|_2^2+ \|b_z(0)\|_2^2\right)^{\frac{3}{2\alpha-3}} +C\,(\|u_0\|_2^2 +\|b_0\|_2^2 + M(t))\right]^{\frac{2\alpha-3}{3}} \label{uzbound}
\end{eqnarray}
and
\begin{eqnarray}
&&
\int_0^t \left(\|\nabla u_z(\tau)\|_2^2 + \|\nabla b_z(\tau)\|_2^2\right)\,d\tau \le F(M(t)) < \infty,  \label{uzbound2}
\end{eqnarray}
where $F(M(t))$ is an explicit function of $M(t)$.
\end{proposition}

\begin{proof}
It is easy to see that $(u,b)$ satisfies
\begin{eqnarray}
&&\hskip-.8in
\|u(t)\|^2_2 + \|b(t)\|_2^2 +  2\;\int_0^t (\|\nabla u (\tau)\|_2^2 + \|\nabla b (\tau)\|_2^2) \,d\tau
\leq \|u_0\|^2_2 + \|b_0\|_2^2.  \label{l2}
\end{eqnarray}
Adding the inner products of $u_z$ with $\pp_z$ of (\ref{DE-1}) and of $b_z$ with $\pp_z$ of (\ref{DE-2}), we obtain, after integration by parts,
\begin{eqnarray*}
&&\hskip-.2in
\frac12\,\frac{d\, (\|u_z\|_2^2+ \|b_z\|_2^2)}{dt} \,
+ \|\nabla u_z\|_2^2 + \|\nabla b_z\|_2^2  \\
&&
=-\,\int
\left[ (u_z \cdot \nabla u) \cdot u_z   - (b_z \cdot \nabla b) \cdot u_z +
(u_z \cdot \nabla b) \cdot b_z -(b_z \cdot \nabla u) \cdot b_z \right] dxdydz\\
&&
\equiv I_1+I_2+I_3+I_4.
\end{eqnarray*}
To bound $I_1$, we integrate by parts and apply H\"{o}lder's inequality to obtain
$$
|I_1| = \left|\int (u_z\cdot \na u_z)\cdot u  \right|
\le C \|\na u_z||_2 \; \|u_z\|_r \; \|u\|_{3\alpha},
$$
where we have omitted $dxdydz$ in the integral for notational convenience and $r$ satisfies
\begin{equation} \label{rrp}
2\le r\le 6,\quad \frac{1}{r}+\frac{1}{3\alpha} = \frac12.
\end{equation}
Applying the Sobolev  inequality
$$
\|u_z\|_r \le C \; \|u_z\|_2^{1- 3(\frac{1}{2} - \frac{1}{r}) }\; \|\na u_z\|_2^{3(\frac{1}{2} - \frac{1}{r}) }
$$
and bounding $\|u\|_{3\alpha}$ by (\ref{3q}), we find
$$
|I_1| \le C\; \|\na u_z||_2^{1+ 3(\frac{1}{2} - \frac{1}{r}) }\; \|u_z\|_2^{1- 3(\frac{1}{2} - \frac{1}{r}) }\; \|u_z\|_{\alpha}^\frac13\; \|\nabla u\|_{2}^\frac23.
$$
By Young's inequality,
\begin{eqnarray*}
|I_1| & \le & \frac14\;\|\na u_z||_2^2 + C\; \|u_z\|_2^2\; \|u_z\|_{\alpha}^q\; \|\na u\|_{2}^{2q}
\end{eqnarray*}
with
\begin{equation}\label{qr}
q= \frac{2}{3-9(\frac12-\frac1r)} = \frac{2}{3(1-\frac1\alpha)}.
\end{equation}
When $\alpha\ge 3$, we have $2q \le 2$ and another application of Young's inequality implies
$$
|I_1| \le  \frac14\; \|\na u_z||_2^2  +  C\; \|u_z\|_2^2 \;\left(\|u_z\|_{\alpha}^{\gamma} + \|\na u\|_{2}^2\right),
$$
where
$$
\gamma\equiv \frac{q}{1-q} = \frac{2}{1-\frac{3}{\alpha}} \quad\mbox{or }\quad \frac{3}{\alpha}+\frac{2}{\gamma}=1.
$$
We now bound $I_2$. By H\"{o}lder's, Sobolev's and Young's inequalities,
\begin{eqnarray*}
|I_2| &\le& C\;\|\na b\|_2 \; \|u_z\|_{\alpha}\; \|b_z\|_{\frac{2\alpha}{\alpha-2}} \\
 &\le& C\;   \|\na b\|_2 \; \|u_z\|_{\alpha}\; \|b_z\|_2^{1-\frac{3}{\alpha}}
 \; \|\na b_z\|_2^{\frac{3}{\alpha}}  \\
  &\le& \frac14\; \|\na b_z\|_2^2 + C\;   \|\na b\|_2^{\frac{2\alpha}{2\alpha-3}} \; \|u_z\|_{\alpha}^{\frac{2\alpha}{2\alpha-3}}\;\| b_z\|_2^{\frac{2\alpha-6}{2\alpha-3}}\\
   &\le& \frac14\; \|\na b_z\|_2^2 + C\;   \left( \|\na b\|_2^2 + \|u_z\|_{\alpha}^\gamma\right)\; \| b_z\|_2^{\frac{2\alpha-6}{2\alpha-3}}
\end{eqnarray*}
where
$$
\gamma =\frac{2\alpha}{\alpha-3} \quad\mbox{or}\quad \frac{3}{\alpha}+\frac{2}{\gamma}=1.
$$
$I_3$ can be bounded exactly as $I_2$. To bound $I_4$, we integrate by parts and apply H\"{o}lder's inequality,
$$
I_4 =-\int
\left[ (b_z \cdot \nabla u) \cdot b_z \right]=  \int \left[ (b_z \cdot \nabla b_z) \cdot u \right]
\le \|\nabla b_z\|_2\; \|b_z\|_r \;\|u\|_{3\alpha},
$$
where $\frac{1}{r}+\frac{1}{3\alpha} = \frac12.$ Following the steps as in the bound of $I_1$, we have
$$
|I_4| \le   \frac14\; \|\na b_z||_2^2 +  C\,\left(\|u_z\|_{\alpha}^{\gamma} + \|\na u\|_{2}^2\right)\;\|b_z\|_2^2.
$$
Combining the estimates for $I_1$, $I_2$, $I_3$ and $I_4$, we find
\begin{eqnarray}
&& \frac{d\, (\|u_z\|_2^2+ \|b_z\|_2^2)}{dt} \,
+ \|\nabla u_z\|_2^2 + \|\nabla b_z\|_2^2
\nonumber\\
&& \qquad \le C \left(\|u_z\|_{\alpha}^{\gamma} + \|\na u\|_{2}^2\right)\; (\|u_z\|_2^2
+ \|b_z\|_2^2) + C  \left( \|\na b\|_2^2 + \|u_z\|_{\alpha}^\gamma\right)\; \| b_z\|_2^{\frac{2\alpha-6}{2\alpha-3}}. \label{fb}
\end{eqnarray}
(\ref{uzbound}) and (\ref{uzbound2}) then follows from (\ref{l2}), (\ref{fb}) and Gronwall's inequality.
\end{proof}

\subsection{Bounds for $\|\na u\|_2$ and $\|\na b\|_2$} This subsection establishes bounds for $\|\na u\|_2$ and $\|\na b\|_2$.

\begin{proposition}
Assume $(u_0, b_0)\in H^3$, $\na \cdot u_0=0$ and $\na\cdot b_0=0$. Let $(u,b)$ be the corresponding solution of the 3D MHD equations (\ref{DE-1}),(\ref{DE-2}),(\ref{DE-3}) and (\ref{DE-4}). Suppose (\ref{uzcon}) holds. Then, for any $t\le T$,
$$
\|\na u(t)\|_2^2 + \|\na b(t)\|_2^2 + \int_0^t (\|\Dd u(\tau)\|_2^2 +\|\Dd b(\tau)\|_2^2)\,d\tau \le G( M(t)) <\infty,
$$
where $G(M(t))$ denotes an explicit function of $M(t)$.
\end{proposition}
\begin{proof}
Adding the inner products of (\ref{DE-1}) with $\Delta u$ and of (\ref{DE-2}) with $\Delta b$ and integrating by parts, we have
\begin{eqnarray}
&&\hskip-.8in
\frac12\,\frac{d}{dt}
\left(\|\na u\|_2^2 + \|\na b\|_2^2\right) + \|\Dd u\|_2^2 + \|\Dd b\|_2^2  \\
&&\hskip-.8in
\qquad = -\int u\cdot \na u \cdot \Dd u + \int b\cdot \na b \cdot \Dd u-\int u\cdot \na b\cdot \Dd b + \int b\cdot \na u \cdot \Dd b.
\end{eqnarray}
By further integrating by parts, we obtain
$$
-\int u\cdot \na u \cdot \Dd u + \int b\cdot \na b \cdot \Dd u-\int u\cdot \na b\cdot \Dd b + \int b\cdot \na u \cdot \Dd b \le \|\na u\|_3^3 + 3 \|\na u\|_3 \;\|\na b\|_3^2.
$$
By (\ref{SL}),
$$
\|\na u\|_3^3 \le C\; \left(\|\na u\|_2^\frac12\;  \|\na_h\na u\|_2^\frac13\; \|\na u_z\|_2^\frac16 \right)^3,
$$
where $\na_h\equiv (\pp_x, \pp_y)$. By Young's inequality,
$$
\|\na u\|_3^3 \le \frac14  \|\na_h\na u\|_2^2 + C\;\|\na u\|_2^3\;\|\na u_z\|_2 \le \frac14  \|\na_h\na u\|_2^2 + C\;\left(\|\na u\|_2^2+\|\na u_z\|_2^2\right)\;\|\na u\|_2^2.
$$
Similarly,
$$
\|\na u\|_3 \;\|\na b\|_3^2 \le \frac14  \|\na_h\na u\|_2^2 + \frac12  \|\na_h\na b\|_2^2 + C\,(\|\na u\|_2^2 + \|\na u_z\|_2^2 + \|\na b_z\|_2^2) \,\|\na b\|_2^2.
$$
Therefore,
\begin{eqnarray*}
\frac{d}{dt}
\left(\|\na u\|_2^2 + \|\na b\|_2^2\right) + \|\Dd u\|_2^2 + \|\Dd b\|_2^2
\le C\; (\|\na u\|_2^2 + \|\na u_z\|_2^2 + \|\na b_z\|_2^2) \,(\|\na u\|_2^2 + \|\na b\|_2^2).
\end{eqnarray*}
Gronwall's inequality coupled with Proposition \ref{ma1} then yields the desired bounds.
\end{proof}

\vspace{.2in}
\section{Criterion in terms of $p_z$}

\setcounter{equation}{0}
\setcounter{theorem}{0}

This section presents the regularity criterion with an assumption on
$p_z$.

\begin{theorem}
Assume the initial data $(u_0, b_0) \in H^1\cap L^4$, $\na\cdot u_0=0$ and $\na\cdot b_0=0$. Let $(u,b)$ be the corresponding solution of the 3D MHD equations (\ref{DE-1}),(\ref{DE-2}),(\ref{DE-3}) and (\ref{DE-4}).   If the pressure $p$ associated with the solution satisfies
\begin{equation}\label{pzcon}
\int_0^T \|p_z(\tau)\|^{\beta}_{\aa} \,d\tau < \infty \quad \mbox{with}\quad \aa \geq \frac{12}{7}\quad\mbox{and}\quad  \frac{3}{\aa} + \frac{2}{\beta} \leq \frac{7}{4}
\end{equation}
for some $T>0$, then $(u,b)$ remains regular on $[0,T]$, namely $(u,b)\in C([0,T]; H^1\cap L^4)$.
\end{theorem}

Since higher-order Sobolev norms of $(u,b)$ can be controlled by its $H^1$-norm (see e.g. \cite{SeTe}), a special consequence of this theorem is that (\ref{pzcon}) yields the global regularity of classical solutions. To prove this theorem, we establish the $L^4$-bound of $(u,b)$ and the desired regularity then follows from the standard Serrin type criteria on the 3D MHD equations \cite{Wu2}.

\begin{proposition} \label{pzp}
Assume the initial data $(u_0, b_0) \in H^1\cap L^4$, $\na\cdot u_0=0$ and $\na\cdot b_0=0$. Let $(u,b)$ be the corresponding solution of the 3D MHD equations (\ref{DE-1}),(\ref{DE-2}),(\ref{DE-3}) and (\ref{DE-4}). If the pressure $p$ satisfies (\ref{pzcon}), then $(u,b)$ obeys the bound
\begin{eqnarray*}
\|w^+\|_4^4 + \|w^-\|_4^4&+& \int_0^t  \left(\|\na |w^+|^2\|_2^2 + \|\na |w^-|^2\|_2^2\right)\, d\tau\\
&+& \, 4\,\int_0^t \int \left(|w^+|^2 \,|\na w^+|^2 + |w^-|^2 \,|\na w^-|^2 \right) \,dxdydz\,d\tau <\infty
\end{eqnarray*}
for any $t\le T$, where
$$
w^\pm = u \pm b.
$$
\end{proposition}

\begin{proof}[Proof of Proposition \ref{pzp}]
We first convert the MHD equations into a symmetric form.
Adding and subtracting (\ref{DE-1}) and (\ref{DE-2}), we find that $w^+$ and $w^-$ satisfy
\begin{eqnarray}
&&\hskip-.8in
\pp_t w^+ + w^-\cdot\nabla w^+ =\Delta w^+ -\nabla p,  \label{w+}  \\
&&\hskip-.8in
\pp_t w^- + w^+\cdot\nabla w^- =\Delta w^- -\nabla p,   \label{w-}   \\
&&\hskip-.8in
\nabla\cdot w^+=0, \quad \nabla \cdot w^-=0. \label{w+div}
\end{eqnarray}
Adding the inner products of (\ref{w+}) with $w^+\,|w^+|^2$ and of (\ref{w-}) with $w^-\,|w^-|^2$ and integrating by parts, we find
\begin{eqnarray}
&& \hspace{-.3in}
\frac14\,\frac{d}{dt} \left( \|w^+\|_4^4 + \|w^-\|_4^4\right) + \frac12 \left(\|\na |w^+|^2\|_2^2 + \|\na |w^-|^2\|_2^2\right) + \int \left(|w^+|^2 \,|\na w^+|^2 + |w^-|^2 \,|\na w^-|^2 \right) \nonumber \\
&& \hspace{-.2in} = J_1 +J_2, \label{mjj}
\end{eqnarray}
where
$$
J_1 =\int p \; w^+\cdot \na |w^+|^2, \quad J_2=\int p \; w^-\cdot \na |w^-|^2.
$$
By H\"{o}lder's inequality,
$$
J_1  \le C \, \|p\|_4\, \|w^+\|_4\, \|\na |w^+|^2\|_2.
$$
We choose $\lambda$ such that
$$
\frac{1}{\alpha}+ \frac{2}{\lambda} = \frac74 \quad\mbox{or}\quad \frac{3\lambda}{2-\lambda\left(1-\frac1{\alpha}\right)}=4.
$$
It then follows from (\ref{3q}) that
$$
\|p\|_4 \le C \, \|p_z\|^\frac13_\alpha\; \|\nabla p\|_\lambda^\frac23.
$$
To further bound $\|\nabla p\|_{\lambda}$, we take the divergence
of (\ref{w+}) to obtain
$$
\Delta p = -\nabla \cdot (w^-\cdot \nabla w^+).
$$
By H\"{o}lder's inequality,
$$
\|\nabla p\|_\lambda \le C\, \|w^-\|_\frac{2\lambda}{2-\lambda}\; \|\nabla w^+\|_2.
$$
Furthermore,  by Sobolev 's inequality,
$$
\|w^-\|_\frac{2\lambda}{2-\lambda}
 = \||w^-|^2\|_{\frac{\lambda}{2-\lambda}}^\frac12 \le C\; \||w^-|^2\|_2^{\frac{3}{\lambda}-\frac74}\; \|\nabla |w^-|^2\|_2^{\frac94-\frac{3}{\lambda}}
  = C\; \|w^-\|_4^{\frac{6}{\lambda}-\frac72} \; \|\nabla|w^-|^2\|_2^{\frac94-\frac{3}{\lambda}},
$$
where we have used the fact that $\frac{12}{7} \le \alpha$ and thus $\lambda \le \frac{12}{7}$. Therefore,
$$
 \|p\|_4 \le C \, \|p_z\|^\frac13_\alpha\; \|\nabla w^+\|^\frac23_2\; \|w^-\|_4^{\frac{4}{\lambda}-\frac73}\; \|\nabla|w^-|^2\|_2^{\frac32-\frac{2}{\lambda}}
$$
and thus
$$
J_1 \le C\;\|p_z\|^\frac13_\alpha\;\|\nabla w^+\|^\frac23_2\; \|w^-\|_4^{\frac{4}{\lambda}-\frac73}\; \|\nabla|w^-|^2\|_2^{\frac32-\frac{2}{\lambda}}\;\|w^+\|_4 \;\|\na |w^+|^2\|_2.
$$
By Young's inequality,
\begin{eqnarray*}
J_1 &\le& \frac{1}{8}\, \|\na |w^+|^2\|_2^2 + \frac{1}{8}\,\|\nabla|w^-|^2\|_2^2
\nonumber\\
&& + \; C\; \|p_z\|_\alpha^{\frac{4\lambda}{3(4-\lambda)}}\;\|\nabla w^+\|_2^{\frac{8\lambda}{3(4-\lambda)}}\;
\|w^-\|_4^{\frac{4(12-7\lambda)}{3(4-\lambda)}}\; \|w^+\|_4^{\frac{4\lambda}{4-\lambda}}.
\end{eqnarray*}
Further applications of Young's inequality imply
\begin{eqnarray*}
&& \|p_z\|_\alpha^{\frac{4\lambda}{3(4-\lambda)}}\;\|\nabla w^+\|_2^{\frac{8\lambda}{3(4-\lambda)}} \le \|p_z\|_\alpha^{\frac{4\lambda}{12-7\lambda}} + \|\nabla w^+\|_2^2, \\
&& \|w^-\|_4^{\frac{4(12-7\lambda)}{3(4-\lambda)}}\; \|w^+\|_4^{\frac{4\lambda}{4-\lambda}} \le \|w^+\|_4^4 + \|w^-\|_4^{\frac{2(12-7\lambda)}{3(2-\lambda)}}.
\end{eqnarray*}
Since $\frac{2(12-7\lambda)}{3(2-\lambda)} <4$, we obtain without loss of generality that \begin{eqnarray}
J_1 &\le& \frac{1}{8}\, \|\na |w^+|^2\|_2^2 + \frac{1}{8}\,\|\nabla|w^-|^2\|_2^2
\nonumber\\
&& + \; C\;\left(\|p_z\|_\alpha^{\frac{4\lambda}{12-7\lambda}} + \|\nabla w^+\|_2^2\right)\,\left(\|w^+\|_4^4 + \|w^-\|_4^4\right). \label{j1b}
\end{eqnarray}
Similarly, we have
\begin{eqnarray}
J_2 &=& \int p \; w^-\cdot \na |w^-|^2 \nonumber\\
&\le& \frac{1}{8}\, \|\na |w^+|^2\|_2^2 + \frac{1}{8}\,\|\nabla|w^-|^2\|_2^2
\nonumber\\
&& + \; C\;\left(\|p_z\|_\alpha^{\frac{4\lambda}{12-7\lambda}} + \|\nabla w^-\|_2^2\right)\,\left(\|w^+\|_4^4 + \|w^-\|_4^4\right). \label{j2b}
\end{eqnarray}
Inserting (\ref{j1b}) and (\ref{j2b}) in (\ref{mjj}) and applying Gronwall's inequality, we obtain the desired result.
\end{proof}

\vskip .4in
\section*{Acknowledgments}
Cao is partially supported by NSF grant DMS 0709228 and a FIU foundation. Wu is partially supported by the AT \& T Foundation at OSU. We thank Professor B. Yuan for careful reading of this manuscript and for discussions.

\vskip .4in
\appendix

\section{}

\renewcommand{\theequation}{\Alph{section}.\arabic{equation}}
\renewcommand{\thetheorem}{\Alph{section}.\arabic{theorem}}
\setcounter{equation}{0}
\setcounter{theorem}{0}

This appendix provides the proofs of the inequalities (\ref{3q}) and (\ref{SL}). For the convenience of future references, we write these inequalities as lemmas.

\begin{lemma} \label{3qin}
Let $\mu$, $\lambda$ and $\gamma$ be three parameters that satisfy
$$
1\le \mu, \;\lambda  <\infty, \quad 1<\frac{1}{\mu}+ \frac{2}{\lambda} \le 4 \quad
\mbox{and}\quad \gamma = \frac{3\lambda}{2-\lambda\left(1-\frac{1}{\mu}\right)}.
$$
Assume $\phi \in H^1({\mathbf R}^3)$, $\phi_x, \; \phi_y\in L^\lambda({\mathbf R}^3)$ and $\phi_z\in L^\mu({\mathbf R}^3)$. Then, there exists a constant $C=C(\mu, \lambda)$ such that
\begin{equation}\label{good0}
\|\phi\|_{\gamma}  \le C\; \|\phi_x\|_\lambda^\frac13\; \|\phi_y\|_\lambda^\frac13 \; \|\phi_z\|^{\frac13}_\mu.
\end{equation}
Especially, when $\lambda=2$, there exists a constant $C=C(\mu)$ such that
\begin{equation}\label{good1}
\|\phi\|_{3\mu}  \le C\; \|\phi_x\|_2^\frac13\; \|\phi_y\|_2^\frac13 \; \|\phi_z\|^{\frac13}_\mu,
\end{equation}
which holds for any $\phi \in H^1({\mathbf R}^3)$ and $\phi_z\in L^\mu({\mathbf R}^3)$ with $1\le \mu<\infty$.
\end{lemma}
\begin{proof}
Clearly,
\begin{eqnarray}
|\phi(x,y,z)|^{1+(1-\frac1\lambda) \gamma} &\leq& C \int_{-\infty}^x |\phi(t,y,z)|^{(1-\frac1\lambda) \gamma} |\pp_t \phi(t,y,z)|\;dt, \label{g111} \\
|\phi(x,y,z)|^{1+(1-\frac1\lambda) \gamma} &\leq& C \int_{-\infty}^y |\phi(x,t,z)|^{(1-\frac1\lambda) \gamma} |\pp_t \phi(x,t,z)|\;dt, \label{g112} \\
|\phi(x,y,z)|^{1+(1-\frac1\mu) \gamma} &\leq& C \int_{-\infty}^z |\phi(x,y,t)|^{(1-\frac1\mu) \gamma} |\pp_t \phi(x,y,t)|\;dt. \label{g113}
\end{eqnarray}
Therefore,
\begin{eqnarray*}
|\phi(x,y,z)|^{\gamma} &\le& C \left[\int_{-\infty}^\infty |\phi(x,y,z)|^{(1-\frac1\lambda) \gamma}|\pp_x \phi(x,y,z)|\;dx\right]^\frac12
\left[\int_{-\infty}^\infty |\phi(x,y,z)|^{(1-\frac1\lambda) \gamma} |\pp_y \phi(x,y,z)|\;dy\right]^\frac12\;  \\
&& \times \left[\int_{-\infty}^\infty |\phi(x,y,z)|^{(1-\frac1\mu) \gamma} |\pp_z \phi(x,y,z)|\;dz\right]^\frac12.
\end{eqnarray*}
Integrating with respect to $x$ and applying H\"{o}lder's inequality, we have
\begin{eqnarray*}
\int_{-\infty}^\infty |\phi(x,y,z)|^{\gamma} \;dx &\le& \left[\int_{-\infty}^\infty |\phi(x,y,z)|^{(1-\frac1\lambda) \gamma} |\pp_x \phi(x,y,z)|\;dx\right]^\frac12 \\
&& \times \left[\int_{-\infty}^\infty\int_{-\infty}^\infty |\phi(x,y,z)|^{(1-\frac1\lambda) \gamma} |\pp_y \phi(x,y,z)|\;dxdy\right]^\frac12 \\
&& \times \left[\int_{-\infty}^\infty\int_{-\infty}^\infty |\phi(x,y,z)|^{(1-\frac1\mu) \gamma} |\pp_z \phi(x,y,z)|\;dxdz\right]^\frac12.
\end{eqnarray*}
Further integration with respect to $y$ and $z$ yields,
\begin{eqnarray*}
\int_{{\mathbf R}^3}|\phi(x,y,z)|^{\gamma} \;dxdydz &\le& \left[\int_{{\mathbf R}^3} |\phi(x,y,z)|^{(1-\frac1\lambda) \gamma} |\pp_x \phi(x,y,z)|\;dxdydz\right]^\frac12 \\
&& \times \left[\int_{{\mathbf R}^3} |\phi(x,y,z)|^{(1-\frac1\lambda) \gamma} |\pp_y \phi(x,y,z)|\;dxdydz\right]^\frac12 \\
&& \times \left[\int_{{\mathbf R}^3} |\phi(x,y,z)|^{(1-\frac1\mu) \gamma} |\pp_z \phi(x,y,z)|\;dxdydz\right]^\frac12.
\end{eqnarray*}
If H\"{o}lder's inequality is applied again, we have
$$
\|\phi\|_{\gamma}^{\gamma} \le C\,\|\phi\|_{\gamma}^{(1-\frac1\lambda)\,\frac{\gamma}{2}}\;\|\pp_x \phi\|_\lambda^\frac12\;
\|\phi\|_{\gamma}^{(1-\frac1\lambda)\,\frac{\gamma}{2}}\;\|\pp_y \phi\|_\lambda^\frac12\; \|\phi\|_{\gamma}^{(1-\frac1\mu)\,\frac{\gamma}{2}}\; \|\pp_z \phi\|_\mu^\frac12,
$$
which leads to (\ref{good0}).
\end{proof}


\begin{lemma}\label{rin}
Let $2\le q\le 6$ and assume $\phi\in H^1({\mathbf R}^3)$. Then, there exists a constant $C=C(q)$ such that
\begin{equation}
\|\phi\|_q \le C\, \|\phi\|_2^{\frac{6-q}{2q}}\; \|\pp_x \phi\|_2^\frac{q-2}{2q}\; \|\pp_y \phi\|_2^\frac{q-2}{2q}\; \|\pp_z \phi\|_2^\frac{q-2}{2q}.
\end{equation}
\end{lemma}
\begin{proof}
This inequality can be obtained by interpolating the trivial inequality $\|\phi\|_2 \le \|\phi\|_2$ and (\ref{good1}) with $\mu=2$, namely
$$
\|\phi\|_6 \le C\; \|\phi_x\|_2^\frac13\; \|\phi_y\|_2^\frac13\; \|\phi_z\|_2^\frac13.
$$
\end{proof}

\end{document}